\documentclass[11pt]{amsart}

\usepackage{amsmath, amssymb, amsfonts, amsthm, bm}

\newcommand{\card}[1]{\lvert#1\rvert}
\DeclareMathOperator{\vol}{vol}
\DeclareMathOperator{\bd}{bd}

\DeclareMathOperator{\diam}{diam}
\newcommand{\norm}[1]{\lVert#1\rVert}
\newcommand{\ipr}[2]{\left\langle #1, #2 \right\rangle}
\newcommand{\numbersystem}[1]{\mathbb{#1}}
\newcommand{\R}{\numbersystem{R}}

\newcommand{\eucl}{\mathbb{E}}
\newcommand{\sph}{\mathbb{S}}
\newcommand{\hyp}{\mathbb{H}}

\newcommand{\abs}[1]{\lvert#1\rvert}

\newcommand{\epsi}{\varepsilon}

\theoremstyle{plain}
\newtheorem{theorem}{Theorem}
\newtheorem{lemma}{Lemma}
\newtheorem*{lmtheorem}{Lawlor-Morgan Theorem}

\begin{document}

\title[Edge-antipodal and subequilateral polytopes]{Upper bounds for edge-antipodal and subequilateral polytopes}
\author{Konrad J. Swanepoel}
\address{Department of Mathematical Sciences,
	University of South Africa, PO Box 392,
	Pretoria 0003, South Africa}
\email{\texttt{swanekj@unisa.ac.za}}
\thanks{This material is based upon work supported by the South African National Research Foundation under Grant number 2053752.}
%\date{\today}

\begin{abstract}
A polytope in a finite-dimensional normed space is \emph{subequilateral} if the length in the norm of each of its edges equals its diameter.
Subequilateral polytopes occur in the study of two unrelated subjects: surface energy minimizing cones and edge-antipodal polytopes.
We show that the number of vertices of a subequilateral polytope in any $d$-dimensional normed space is bounded above by $(\frac{d}{2}+1)^d$ for any $d\geq 2$.
The same upper bound then follows for the number of vertices of the edge-antipodal polytopes introduced by I.~Talata (Period.\ Math.\ Hungar.\ \textbf{38} (1999), 231--246).
This is a constructive improvement to the result of A.~P\'or (to appear) that for each dimension $d$ there exists an upper bound $f(d)$ for the number of vertices of an edge-antipodal $d$-polytopes.
We also show that in $d$-dimensional Euclidean space the only subequilateral polytopes are equilateral simplices.
\end{abstract}

\maketitle

\section{Notation}
Denote the $d$-dimensional real linear space by $\R^d$, a norm on $\R^d$ by $\norm{\cdot}$,
its unit ball by $B$, and the ball with centre $x$ and radius $r$ by $B(x,r)$.
Denote the diameter of a set $C\subseteq \R^d$ by $\diam(C)$, and (if it is measurable) its volume (or $d$-dimensional Lebesgue measure) by $\vol(C)$.
The \emph{dual norm} $\norm{\cdot}^\ast$ is defined by $\norm{x}^\ast:=\sup\{\ipr{x}{y} : \norm{y}\leq 1\}$, where $\ipr{\cdot}{\cdot}$ is the inner product on $\R^d$.
Denote the number of elements of a finite set $S$ by $\card{S}$.
The \emph{difference body} of a set $S\subseteq\R^d$ is $S-S:=\{x-y:x,y\in S\}$.
A \emph{polytope} is the convex hull of finitely many points in some $\R^d$.
A \emph{$d$-polytope} is a polytope of dimension $d$.
A \emph{convex body} $C$ is a compact convex subset of $\R^d$ with nonempty interior.
The boundary of $C$ is denoted by $\bd C$.
Given any convex body $C$ we define the \emph{relative norm $\norm{\cdot}_C$ determined by $C$} to be the norm with unit ball $C-C$, or equivalently,
\[ \norm{x}_C := \sup\{\lambda>0: a+\lambda x\in C \text{ for some } a\in C\}.\]
See \cite{Grunbaum, Ba, MR97f:52001} for background on polytopes, convexity, and finite-dimensional normed spaces.

\section{Introduction}

\subsection{Antipodal and edge-antipodal polytopes}
A $d$-polytope $P$ is \emph{antipodal} if for any two vertices $x$ and $y$ of $P$ there exist two parallel hyperplanes, one through $x$ and one through $y$, such that $P$ is contained in the closed slab bounded by the two hyperplanes.
Klee \cite{Klee} posed the problem of finding an upper bound for the number of vertices of an antipodal $d$-polytope in terms of $d$.
Danzer and Gr\"unbaum \cite{MR25:1488} proved the sharp upper bound of $2^d$.
See \cite{MS} for a recent survey.

A $d$-polytope $P$ is \emph{edge-antipodal} if for any two vertices $x$ and $y$ joined by an edge there exist two parallel hyperplanes, one through $x$ and one through $y$, such that $P$ is contained in the closed slab bounded by the two hyperplanes.
This notion was introduced by Talata \cite{Talata}, who conjectured that the number of vertices of an edge-antipodal $3$-polytope is bounded above by a constant.
Csik\'os \cite{Csikos} proved an upper bound of $12$, and
K.~Bezdek, Bisztriczky and B\"or\"oczky \cite{BBB} gave the sharp upper bound of $8$.
P\'or \cite{Por} proved that the number of vertices of an edge-antipodal $d$-polytope is bounded above by a function of $d$.
However, his proof is existential, with no information on the size of the upper bound.
Our main result is an explicit bound.

\begin{theorem}\label{th1}
Let $d\geq 2$.
Then the number of vertices of an edge-antipodal $d$-polytope is bounded above by $(\frac{d}{2}+1)^d$.
\end{theorem}

In the plane, an edge-antipodal polytope is clearly antipodal, and in this case the above theorem is sharp.
The bound given is not sharp for $d\geq 3$ (since the bound in Theorem~\ref{th2} below is not sharp).
In \cite{BBB} it is stated without proof that all edge-antipodal $3$-polytopes are antipodal.
On the other hand, Talata has an example of an edge-antipodal $d$-polytope that is not antipodal for each $d\geq 4$ (see \cite{Csikos} and Section~\ref{s4} below).
Most likely the largest number of vertices of an edge-antipodal $d$-polytope has an upper bound exponential in $d$, perhaps even $2^d$.
We also mention the paper by Bisztriczky and B{\"o}r{\"o}czky \cite{BB} discussing edge-antipodal $3$-polytopes.

Theorem~\ref{th1} is proved by considering a metric relative of edge-antipodal polytopes, discussed next.

\subsection{Equilateral and subequilateral polytopes}
A polytope $P$ is \emph{equilateral} with respect to a norm $\norm{\cdot}$ on $\R^d$ if its vertex set is an \emph{equidistant set}, i.e., the distance between any two vertices is a constant.
This notion was first considered by Petty \cite{MR43:1051}, who showed that equilateral polytopes are antipodal, hence have at most $2^d$ vertices.
We now introduce the following natural weakening of this notion, analogous to the weakening from antipodal to edge-antipodal.
We say that a $d$-polytope $P$ is \emph{subequilateral} with respect to a norm $\norm{\cdot}$ on $\R^d$ if the length of each of its edges equals its diameter.

Although not explicitly given a name, the vertex sets of subequilateral polytopes appear in the study of surface energy minimizing cones by Lawlor and Morgan \cite{MR95i:58051}; see Section~\ref{s4} for a discussion.

It is well-known and easy to prove that an edge-antipodal polytope $P$ is subequilateral with diameter $1$ in the relative norm $\norm{\cdot}_P$ determined by $P$ \cite{Talata, Csikos}.
It is also easy to see that any subequilateral polytope is edge-antipodal.
In order to prove Theorem~\ref{th1} it is therefore sufficient to bound the number of vertices of a subequilateral $d$-polytope.

\begin{theorem}\label{th2}
Let $d\geq 2$.
Then the number of vertices of a subequilateral $d$-polytope with respect to some norm $\norm{\cdot}$ is bounded above by $(\frac{d}{2}+1)^d$.
\end{theorem}
The proof is in Section~\ref{s2}.
In two-dimensional normed spaces subequilateral polytopes are always equilateral.
Therefore, the above theorem is sharp for $d=2$.
By analyzing equality in the proof of Theorem~\ref{th2}, it can be seen that the bound is not sharp for $d\geq 3$.
Since edge-antipodal $3$-polytopes have at most $8$ vertices, with equality only for parallelepipeds \cite{BBB}, it follows that subequilateral $3$-polytopes with respect to any norm has size at most $8$, with equality only if the unit ball of the norm is a parallelepiped homothetic to the polytope.

We finally mention that in Euclidean $d$-space $\eucl^d$ the only subequilateral polytopes are equilateral simplices, and give a proof.
In the proof we have to consider subequilateral polytopes in spherical spaces, making it possible to formulate a more general theorem for spaces of constant curvature.
Note that if we restrict ourselves to a hemisphere of the $d$-sphere $\sph^d$ in $\eucl^{d+1}$, the notion of a polytope can be defined without ambiguity.
The definition of a subequilateral polytope then still makes sense in in a hemisphere of $\sph^d$, as well as in hyperbolic $d$-space $\hyp^d$.
\begin{theorem}\label{th3}
Let $P$ be a subequilateral $d$-polytope in either $\eucl^d$, $\hyp^d$, or a hemisphere of $\sph^d$.
Then $P$ is an equilateral $d$-simplex.
\end{theorem}
\begin{proof}
The proof is by induction on $d\geq 1$, with $d=1$ trivial and $d=2$ easy.
Suppose now $d\geq 3$.
Let $P$ be a subequilateral $d$-polytope in any of the three spaces.
By induction all facets of $P$ are equilateral simplices.
In particular, $P$ is simplicial.
Since $d\geq 3$, it is sufficient to show that $P$ is simple (see section~4.5 and exercise~4.8.11 of \cite{Grunbaum}).

Consider any vertex $v$ with neighbours $v_1,\dots,v_k$, $k\geq d$.
Then $v_1,\dots,v_k$ are contained in an open hemisphere $S$ of the $(d-1)$-sphere of radius $\diam(P)$ and centre $v$.
(This sphere will be isometric to some sphere in $\eucl^{d}$, not necessarily of radius $\diam(P)$.)

Consider the $(d-1)$-polytope $P'$ in $S$ generated by $v_1,\dots,v_k$ and any facet of $P'$ with vertex set $F\subset\{v_1,\dots,v_k\}$.
There exists a great sphere $C$ of $S$ passing through $F$ with $P'$ in one of the closed hemispheres determined by $C$.
It follows that the hyperplane $H$ generated by $C$ and $v$ passes through $F\cup\{v\}$, and $P$ is contained in one of the closed half spaces bounded by $H$.
Therefore, $F\cup\{v\}$ is the vertex set of a facet of $P$.

Similarly, it follows that for any vertex set $F$ of a facet of $P$ containing $v$, $F\setminus\{v\}$ is the vertex set of a facet of $P'$.
Therefore, any edge $v_iv_j$ of $P'$ is an edge of $P$, hence of length the diameter of $P$.
It follows that the distance between $v_i$ and $v_j$ in $H$ is the diameter of $P'$ as measured in $H$.
This shows that $P'$ is subequilateral in $H$, and so by induction is an equilateral $(d-1)$-simplex.
Therefore, $k=d$, giving that $P$ is a simple polytope, which finishes the proof.
\end{proof}

\section{A measure of non-equidistance}\label{s2}
The key to the proof  of Theorem~\ref{th2} is a lower bound for the distance between two nonadjacent vertices of a subequilateral polytope.
For any finite set of points $V$ we define
\[\lambda(V;\norm{\cdot})=\diam(V)/\min_{x,y\in V, x\neq y}\norm{x-y}.\]
Since $\lambda(V;\norm{\cdot})\geq 1$, with equality if and only if $V$ is equidistant in the norm $\norm{\cdot}$, this functional measures how far $V$ is from being equidistant.
The next lemma generalizes the theorem of Petty \cite{MR43:1051} and Soltan \cite{MR52:4127} that the number of points in an equidistant set is bounded above by $2^d$.
In \cite{MR93d:52009} a proof of the $2^d$-upper bound was given using the isodiametric inequality for finite-dimensional normed spaces due to Busemann (equation~(2.2) on p.~241 of \cite{B}; see also Mel'nikov \cite{MR27:6191}).
However, since the isodiametric inequality has a quick proof using the Brunn-Minkowski inequality \cite{BZ}, it is not surprising that the latter inequality occurs in the following proof.
\begin{lemma}\label{l2}
Let $V$ be a finite set in a $d$-dimensional normed space.
Then $\card{V}\leq(\lambda(V;\norm{\cdot})+1)^d$.
\end{lemma}
\begin{proof}
Let $\lambda=\lambda(V;\norm{\cdot})$.
By scaling we may assume that $\diam(V)=\lambda$.
Then $\norm{x-y}\geq 1$ for all $x,y\in V$, $x\neq y$, hence the balls $B(v,1/2)$, $v\in V$, have disjoint interiors.
Define $C=\bigcup_{v\in V}B(v,1/2)$.
Then $\vol(C)=\card{V}(1/2)^d\vol(B)$ and $\diam(C)\leq 1+\lambda$.
By the Brunn-Minkowski inequality \cite{BZ} we obtain $\vol(C-C)^{1/d}\geq\vol(C)^{1/d}+\vol(-C)^{1/d}$.
Noting that $C-C\subseteq(1+\lambda)B$, the result follows.
%By the isodiametric inequality for finite-dimensional normed spaces, due to Busemann \cite[(2.2) on p.~241]{B} and Mel'nikov \cite{MR27:6191} (see also \cite[p.~93]{BZ}), balls have the largest volume among measurable sets of a fixed diameter: $\vol(C)\leq\vol(\frac{1+\lambda}{2}B)$.
%Therefore, $\card{S}(\lambda/2)^d\vol(B)\leq(\frac{1+\lambda}{2})^d\vol(B)$, and $\card{S}\leq(\frac{1}{\lambda}+1)^d$.
\end{proof}

In order to find an upper bound on the number of vertices of a subequilateral polytope with vertex set $V$, it remains to bound $\lambda(V;\norm{\cdot})$ from above.

\begin{lemma}\label{l3}
Let $d\geq 2$ and let $V$ be the vertex set of a subequilateral $d$-polytope.
Then $\lambda(V;\norm{\cdot})\leq d/2$.
\end{lemma}

\begin{proof}%[Proof of Lemma~\ref{l3}]
Let $P$ be a subequilateral $d$-polytope of diameter $1$, and let $V$ be its vertex set.
We have to show that $\norm{x-y}\geq2/d$ for any distinct $x,y\in V$.
Since this follows from the definition if $xy$ is an edge of $P$,
%If $xy$ is not an edge, but still on the boundary of $P$, then $d\geq 3$ and we may use induction on $d$, since the vertices of the face of $P$ containing $xy$ is subequilateral in a space of smaller dimension.
we assume without loss that $xy$ is not an edge of $P$.
Then $xy$ intersects the convex hull $P'$ of $V\setminus\{x,y\}$ in a (possibly degenerate) segment, say $x'y'$, with $x$, $x'$, $y'$, $y$ in this order on $xy$.
Let $F_x$ and $F_y$ be facets of $P'$ containing $x'$ and $y'$, respectively.

We show that $\norm{x-x'}\geq1/d$.
For each vertex $z$ of $F_x$, $xz$ is an edge of $P$, hence $\norm{x-z}=1$.
By Carath\'eodory's theorem \cite[(2.2)]{Ba}, there exist $d$ vertices $z_1,\dots,z_d$ of the $(d-1)$-polytope $F_x$ and real numbers $\lambda_1,\dots,\lambda_d$ such that
\[ x'=\sum_{i=1}^d\lambda_i z_i,\quad \lambda_i\geq 0,\quad \sum_{i=1}^d\lambda_i=1.\]
Suppose without loss that $\lambda_d=\max_i\lambda_i$.
Then $\lambda_d\geq 1/d$.
By the triangle inequality we obtain
\begin{align*}
\norm{x'-z_d} & = \norm{\sum_{i=1}^{d-1}\lambda_i(z_i-z_d)}\leq \sum_{i=1}^{d-1}\lambda_i\norm{z_i-z_d}\\
&\leq \sum_{i=1}^{d-1}\lambda_i=1-\lambda_d\leq 1-\frac{1}{d},
\end{align*}
and
\begin{align*}
\norm{x-x'} & \geq\norm{x-z_d}-\norm{x'-z_d}\\
& \geq 1-(1-\frac{1}{d})=\frac{1}{d}.
%1=\norm{x-z_d} &\leq \norm{x-x'}+\norm{x'-z_d}\\
%&\leq \norm{x-x'}+1-\frac{1}{d},
\end{align*}
%and $\norm{x-x'}\geq 1/d$.
Similarly, $\norm{y-y'}\geq1/d$, and we obtain $\norm{x-y}\geq 2/d$.
%
%It follows that all nonzero distances between the vertices of a subequilateral polytope of diameter $1$ are in the interval $[2/d,1]$.
%The bound on the size of $S$ now follows from the following standard lemma.
\end{proof}
Lemmas~\ref{l2} and \ref{l3} now imply Theorem~\ref{th2}.\qed

\section{Concluding remarks}\label{s4}
\subsection{Sharpness of Lemma~\ref{l3}}
The following example shows that Lemma~\ref{l3} cannot be improved in general.
Consider the subspace $X=\{(x_1,\dots,x_{d+1}):\sum_{i=1}^dx_i=0\}$ of $\R^{d+1}$ with the $\ell_1$ norm $\norm{(x_1,\dots,x_{d+1})}_1:=\sum_{i=1}^{d+1}\abs{x_i}$.
Let the standard unit vector basis of $\R^{d+1}$ be $e_1,\dots,e_{d+1}$.
Let $c=\sum_{i=1}^d e_i$.
Then $V=\{de_i-c:i=1,\dots,d\}\cup\{\pm 2e_{d+1}\}$ is the vertex set of a $d$-polytope $P$ in $X$, with all intervertex distances equal to $2d$, except for the distance between $\pm 2e_{d+1}$, which is $4$.
It follows that $P$ is subequilateral and $\lambda(V;\norm{\cdot})=d/2$.

However, the above polytope $P$ is in fact antipodal, and so it is equilateral in $\norm{\cdot}_P$, which gives $\lambda(V;\norm{\cdot}_P)=1$.
It is easy to see that for any polytope $P$ subequilateral with respect to some norm $\norm{\cdot}$, and with vertex set $V$, we have $\lambda(V,\norm{\cdot})\leq\lambda(V,\norm{\cdot}_P)$.
One may therefore hope that for the norm $\norm{\cdot}_P$ the upper bound in Lemma~\ref{l3} may be improved, thus giving a better bound in Theorem~\ref{th1}.
The following example shows that any such improved upper bound will still have to be at least  $(d-1)/2$, indicating that essentially new ideas will be needed to improve the upper bounds in Theorems~\ref{th1} and \ref{th2}.

We consider Talata's example \cite{Csikos} of an edge-antipodal polytope that is not antipodal.
Let $d\geq 4$, $e_1,\dots,e_{d}$ be the standard basis of $\R^d$, $p=\frac{2}{d-1}\sum_{i=1}^{d-1}e_i$, and $\lambda=(d-1)/2-\epsi>1$ for some small $\epsi>0$.
Then the polytope $P$ with vertex set $V=\{o,e_1,\dots,e_d,p,e_d+\lambda p\}$ is edge-antipodal but not antipodal.
In fact, $\diam(V)\leq 1$ by definition of $\norm{\cdot}_P$, and since $\norm{e_d-o}_P=1$ and $\norm{p-o}_P=1/\lambda$, we obtain $\lambda(V,\norm{\cdot}_P)\geq\lambda$, which is arbitrarily close to $(d-1)/2$.

\subsection{Subequilateral polytopes in the work of Lawlor and Morgan}
Define the \emph{$\norm{\cdot}$-energy} of a hypersurface $S$ in $\R^d$ to be $\norm{S}:=\int_S\norm{n(x)}dx$, where $n(x)$ is the Euclidean unit normal at $x\in S$.
In \cite{MR95i:58051} a sufficient condition is given to obtain an energy minimizing hypersurface partitioning a convex body.
We restate a special case of the ``General Norms Theorem I'' in \cite[pp.~66--67]{MR95i:58051} in terms of subequilateral polytopes.
(In the notation of \cite{MR95i:58051} we take all the norms $\Phi_{ij}$ to be the same.
Then the points $p_1,\dots,p_m$ in the hypothesis form an equidistant set with respect to the dual norm.
The weakening of the hypothesis in the last sentence of the General Norms Theorem I is easily seen to be equivalent to the requirement that $p_1,\dots,p_m$ is the vertex set of a subequilateral polytope.)
We refer to \cite{MR95i:58051} for the simple and enlightening proof using the divergence theorem.

\begin{lmtheorem}
Let $\norm{\cdot}$ be a norm on $\R^n$, and let $p_1,\dots,p_m\in\R^n$ be the vertex set of a subequilateral polytope of $\norm{\cdot}$-diameter $1$.
Let $\Sigma=\bigcup H_{ij}\subset C$ be a hypersurface which partitions some convex body $C$ into regions $R_1,\dots,R_m$ with $R_i$ and $R_j$ separated by a piece $H_{ij}$ of a hyperplane such that the parallel hyperplane passing through $p_i-p_j$ supports the unit ball $B$ at $p_i-p_j$.

Then for any hypersurface $M=\bigcup M_{ij}$ which also separates the $R_i\cap\bd C$ from each other in $C$, with the regions touching $R_i\cap\bd C$ and $R_j\cap\bd C$ facing each other across $M_{ij}$, we have $\norm{\Sigma}^\ast\leq\norm{M}^\ast$, i.e.\ $\Sigma$ minimizes $\norm{\cdot}^\ast$-energy, where $\norm{\cdot}^\ast$ is the norm dual to $\norm{\cdot}$.
\end{lmtheorem}

%\bibliography{edge-antipodal}

\end{document}